\documentclass[11pt]{article}
\usepackage{amsfonts,euscript,theorem}
\newcommand{\no}{\noindent}
\newcommand{\bc}{\begin{center}}
\newcommand{\ec}{\end{center}}
\newcommand{\be}{\begin{equation}}
\newcommand{\ee}{\end{equation}}
\newcommand{\bea}{\begin{eqnarray*}}
\newcommand{\eea}{\end{eqnarray*}}
\newcommand{\bean}{\begin{eqnarray}}
\newcommand{\eean}{\end{eqnarray}}
\newcommand{\ran}{\ensuremath{\rangle}}
\newcommand{\lan}{\ensuremath{\langle}}
\newcommand{\rar}{\ensuremath{\rightarrow}}

\newcommand{\pd}{\ensuremath{\partial}}
\newcommand{\pdp}{\ensuremath{\partial}^p}
\newcommand{\pdx}{\ensuremath{\partial}_x}
\newcommand{\h}[1]{\ensuremath{\hbox{ #1 }}}

\newcommand{\w}[1]{\ensuremath{\widetilde{#1}}}

\newcommand{\de}{\ensuremath{\delta}}
\newcommand{\te}{\ensuremath{\theta}}
\newcommand{\ga}{\ensuremath{\Gamma}}
\newcommand{\Om}{\ensuremath{\Omega}}
\newcommand{\si}{\ensuremath{\sigma}}

\newcommand{\ps}{\ensuremath{\psi}}
\newcommand{\ve}{\ensuremath{\varepsilon}}
\newcommand{\vp}{\ensuremath{\varphi}}
\newcommand{\pe}{\ensuremath{\vp_{\ve}}}

\newcommand{\N}{\ensuremath{\mathbb{N}}}
\newcommand{\NN}{\ensuremath{\mathbb{N}_0}}
\newcommand{\R}{\ensuremath{\mathbb{R}}}
\newcommand{\C}{\ensuremath{\mathbb{C}}}
\newcommand{\A}{\ensuremath{A_0(\R)}}
\newcommand{\Aq}{\ensuremath{A_q(\R)}}
\newcommand{\G}{\ensuremath{{\EuScript{G}}}}
\newcommand{\GR}{\ensuremath{{\EuScript{G}}(\R)}}
\newcommand{\D}{\ensuremath{{\EuScript{D}}(\R)}}
\newcommand{\DD}{\ensuremath{{\EuScript{D}}'(\R)}}

\theorembodyfont{\itshape}
\newtheorem{Th1}{Theorem}
\newtheorem{Cor1}{ Corollary}
\newtheorem{Th2}[Th1]{Theorem}

\begin{document}
\setlength{\baselineskip}{18pt} \vspace*{2mm}\bc {\bf {RESULTS ON
MODELLING AND PRODUCTS \\ OF SINGULARITIES IN COLOMBEAU ALGEBRA
\boldmath $\GR$ \unboldmath }}

{\textsc{Blagovest Damyanov}}\footnote{ E-mail: bdamyanov@mail.bg}
\\ \textit{Bulgarian Acad. of Sci.,  INRNE - Theor.
Math. Physics Dept. \\ 72 Tzarigradsko shosse, 1784 Sofia, Bulgaria}
\ec

\vspace*{2mm} \setlength{\baselineskip}{13pt} {\small Modelling of
singularities given by discontinuous functions or distributions by
means of generalized functions has proved useful in many problems
posed by physical phenomena.  We introduce in a systematic way
generalized functions of Colombeau that model such singularities.
Moreover, we evaluate some products of singularity-modelling
generalized functions whenever the result admits an associated
distribution.} \setlength{\baselineskip}{18pt} 


\vspace*{3mm}

\no \textbf{1 \ INTRODUCTION}

\vspace*{2mm} The generalized functions of Colombeau \cite{col84}
have  become a powerful tool for treating differential equations
with singular coefficients and data \,as well as singular products
of Schwartz distributions. The flexibility of Colombeau theory
allows to model such singularities by means of appropriately chosen
generalized functions, treat them in this framework and --- via the
association process --- bring down the obtained results to
distributional level.

In particular, generalized models of Heaviside step-function \te \
have proved useful in solving problems that arise in Mathematical
Physics \cite{col87}. Other examples involving \te- \,and \de-type
singularities that describe realistic physical phenomena are jump
conditions in hyperbolic systems leading to travelling \de-waves
solutions \cite{col89}, the so-called controlled hybrid systems
\cite{HO}, geodesics for impulsive gravitational waves  \cite{St}.
A detailed presentation of results and citations on this topic can
be found in \cite{mob} and \cite{gct}.

Motivated by such works, we have introduced in a \textit{unified way} generalized functions of Colombeau
 that model singularities of certain type and have additional properties  \cite{dam}. Following the same idea, we model
 in this paper singularities  given by functions that have discontinuities of first order in a point on the real line \R.
 Moreover, we evaluate products of singularity-modelling generalized functions whenever the result admits an associated distribution.

\newpage
\vspace*{2mm}
\no \textbf{2 \ NOTATION AND DEFINITIONS}
\vspace*{3mm}

\no We recall first the basic definitions of Colombeau algebra $\G(\R)$.

\vspace*{1mm} \textit{Notation 1.} Let \N \ denote the natural
numbers, $\NN = \N \cup \{0\}$, and \,$\de _{ij} = \{\,1$ \,if
\,$i=j,\, = 0$ \,if $i\ne j$\,\}, \,for \,$i, j\in\NN$. Then we
put for arbitrary $q \in \NN$\,:
\[\Aq = \{ \vp(x) \in \D: \int_{\R} x^{j}\,\vp (x)\,dx  = \de _{0j}, \ j = 0, 1,...,q \},\]
where $\D$ is the space of infinitely-differentiable functions
with compact support. For $\vp \in\A$ and $\ve
> 0$, we will use the following notation throughout the paper: \
$\pe = \ve ^{-1}\vp(\ve ^{-1}x)$ \ and \ $s \equiv s(\vp) := sup\
\{|x|: \vp (x)\ne 0 )\}$. Then clearly\ $s(\pe) =\ve s(\vp)$, and
denoting $\si \equiv \si(\vp, \ve) := s(\pe) > 0$, we have $\si :=
\ve s = O(\ve)$, as $\ve \rightarrow 0$, for each $\vp \in\A$. \
Finally, the shorthand notation $\pdx = d/dx$ \,will be used in
the one-dimensional case too.

\vspace*{1mm}
\textsc{Definition 1.} Let $\EuScript{E}\,[\R]$ be the algebra of functions $F(\vp , x): \A \times \R  \rar
\C$ that are infinitely differentiable for  fixed `parameter' \vp. Then the generalized functions of Colombeau are elements of the quotient algebra
\ \ $\G \equiv\GR = \EuScript{E}_{\mathrm{M}} [\R]\,/ \ \EuScript{I}\,[\R]$. \ Here $\EuScript{E}_{\mathrm{M}}[\R]$
\,is the subalgebra of `moderate' functions such that for each compact subset $K$ of \R \,and $p \in\N$ \,there is a $q\in \N$ \,such that, for each
 $\vp \in \Aq$,  \  $\sup_{x \in K}\,|\pdp \,F(\pe, x)\,| = O(\ve^{-q}), \hbox{ as}  \ \ve \rar 0_+$.\ The ideal $\EuScript{I}\,[\R]$
 of $\EuScript{E}_{\mathrm{M}}[\R]$ consists of all functions such that for each compact $K \subset\R $ \ and any $p\in \N$ \,there is a $q\in\N$
 \ such that, for every $r \geq q$ and $\vp \in A_r(\R)$, \ $\sup_{x \in K}\,| \pdp \,F(\pe, x)\,| = O(\ve^{r-q}), \hbox{ as} \ \ve \rar 0_+$.

The algebra \GR \ contains the distributions on \R, canonically embedded as a \C-vector subspace \,by the map  \vspace*{-1mm}
\[ i : \DD \rar \,\G : u \mapsto \w{u} = \{\,\w{u}(\vp, x) := (u * \check{\vp})(x)| \,\vp \in\Aq \,\}, \  \h{where} \ \check{\vp}(x) = \vp(-x).  \]

The equality of generalized functions in  \G \ is very strict and a weaker form of equality in the sense of {\em{association}} is introduced,
which plays a fundamental role in Colombeau theory.

\vspace*{1mm}
\no \textsc{Definition} 2 (a) Two generalized functions $F, G \in\GR$ \,are said to be `associated', denoted $F \approx G$, if for some
representatives $F(\vp _\ve, x), G(\vp _\ve, x)$ and arbitrary $\ps (x)\in\D$ \,there is a $q\in\NN$, such that for any $\vp (x)\in\Aq,
\ \lim_{{}\ve \rar 0_+} \int_{\R} [F(\vp _\ve , x) - G(\vp _\ve, x)] \ps (x)\,dx = 0.$

(b) A generalized function $F\in\GR$ \ is said to be `associated' with a distribution $ u\in \DD$, denoted $F \approx u$, \,if for some
representative $F(\vp _\ve, x)$, \,and arbitrary $\ps (x)\in\D$ \,there is a $q\in\NN$, such that for any $\vp (x)\in\Aq, \ \lim_{{}\ve \rar 0_+}
 \int_{\R} F(\vp _\ve , x) \ps (x)\,dx = \lan u, \ps\ran.$

These definitions are independent of the representatives chosen, and the association is a faithful generalization of the equality of distributions.
The following relations hold in \G\,:
\[ F \approx u \quad  \& \quad F_1 \approx u_1 \  \Longrightarrow \ F + F_1 \approx u + u_1, \quad \pd F \approx \pd u. \]

\textbf{Remark.} Note that the relation $F\approx u$ is asymmetric in the sense that the terms cannot be moved over the $\approx$-sign: on the r.h.s.
of it there stands a distribution. Of course, the equivalent relation  $F\approx \w{u}$ \ in \G \ \,is symmetric (and can be written as
$F - \w{u}\approx 0$ as well). We prefer however the first, simpler and suggesting, notation for the associated distribution.

\vspace*{1mm}
Below we shall need also the following.

\vspace*{1mm}
\textit{Notation 2.} If $a\in \R, \,a>-1$, introduce the locally-integrable functions that are `normed' powers of the variable $x \in \R$
\,supported in one of the real semiaxes\,:
\[ \nu_+^{\,a} \equiv \nu_+^{\,a}(x) =  \{\frac{x^{\,a}}{\ga(a+1)} \ \h{if} x > 0 , \ \ = 0 \ \h{if} x<0\}. \]
\[\nu_-^{\,a} \equiv \nu_-^{\,a}(x)  =  \{\frac{(-x)^{\,a}}{\ga(a+1)} \ \h{if} x < 0,  \  \  = 0 \ \h{if} x >0 \}.  \]
In particular, for $a = p \in \NN$, we have
\[ \nu_+^{\,p} =  \{\ \frac{x^p}{p!}, \ x > 0; \ \ = 0, \  x<0 \}, \qquad \nu_-^{\,p} =  \{\ \frac{ (-x)^p }{p!}, \ x < 0;  \ \ = 0, \ x>0 \}. \]
We can also define their `even' and `odd' compositions\,: \  $|\nu|^{\,p} = \nu_+^{\,p} + \nu_-^{\,p}, \ |\nu|^{\,p}\,\mathrm{sgn}\,x =
\nu_+^{\,p} - \nu_-^{\,p}.$
If $p=0$, we come to the Heaviside step-function \,$\te = \{1, \ x < 0;  \ \ = 0, \ x>0 \}$.

All these functions have singularities at $x=0$ and only considered as distributions, their derivatives exist and satisfy the relations
\[ \pdx \,\nu_{\pm}^{\,p+1} = \pm \ \nu_{\pm}^{\,p}, \quad \pdx \,|\nu|^{\,p} =  |\nu|^{\,p-1}\,\mathrm{sgn}\,x,
\quad \pdx \,|\nu|^{\,p}\,\mathrm{sgn}\,x  =  |\nu|^{\,p-1}. \]
Note that, due to the norming, no number coefficients are present here, which makes the calculations with these distributions easier.

The definition of  $\nu_{\pm}^{\,a}$ \,can be extended for any $a\in\Om := \R\backslash - \N \}$, by setting
\[ \nu_{+}^{\,a} = \pd^r \,\nu_{+}^{\,a +r }(x), \qquad  \nu_{-}^{\,a} = (-1)^r \,\pd^r \,\nu_{-}^{\,a +r }(x), \ \h{ where} \ r\in \NN
\ \h{is such that} \ a+r > -1.\]

\vspace*{1mm}
\textit{Notation 3.}\ Let $C_d^k (\R\setminus\!\{0\})$ be the class of $k$-times differentiable functions on $\R\setminus\!\{0\}$ \
for some $k\in\N$, \,such that each function $f(x)$ and its derivatives $f^{(i)}, i= 1,...,k$, have discontinuities of first order at the
 point $x=0$, i.e. for each $i \leq k$, the values $f^{(i)}(0_{+})$ and $f^{(i)}(0_{-})$ exist \,but in general differ from each other. Then we
 denote the mean values and the jumps at $x=0$ of $f^{(i)}(x)$ \,by\,:
\[m_i \equiv m_i (f) = \frac{1}{2}\,\left[f^{(i)} (0_+) + f^{(i)} (0_-)\right], \quad h_i  \equiv h_i (f) = f^{(i)} (0_+) - f^{(i)} (0_-),
\ i=0,...,k.\]

\vspace*{5mm}
\no\textbf{3 \ MODELLING OF SINGULARITIES IN} {\boldmath{\GR}}\unboldmath
\vspace*{3mm}

Consider first generalized functions that model the \de-type singularity in the sense of association, i.e. being associated with the \de-function.
Since there is an abundant variety of such functions (together with the canonical imbedding $\w{\de}$ in \G \ of the distribution \de), we can put
on the generalized functions in question an additional requirement. So define, following \cite[\S 10]{mob}, a generalized function $D \in \G$
 with the properties:
\be D \approx \de, \qquad D^2 \approx \de. \label{DD} \ee

Let $\vp \in\A$, $s \equiv s(\vp)$, and $\si = s(\pe) = \ve s$ be as in Notation~1, and $D \in \G$ be the class $[ \vp \mapsto D( s(\vp), x)]$. We
 then specify that $D(s, x) = f(x) + \lambda_{s}\ g(x)$, where $f, g \in \D$ are real-valued, symmetric, with disjoint support, and satisfying:
\[ \int_{\R}f(x) dx = 1, \quad \int_{\R}g(x) dx = 0, \ \h{and} \ \lambda_{s}^2 = \frac{s - \int f^2(x) dx}{\int g^2(x) dx}.\]
Now, it is not difficult to check that, for each $\vp \in \A$, the representative $D(s, x)$ of the generalized function $D$ satisfies the conditions:
\be D(., x)\in \D, \quad D(., -x ) = D(., x), \quad \frac{1}{s} \int_{\R} D^{\,2}(s, x)\ dx \ = \ \int_{\R} D(s, x) dx \ = \ 1, \label{D-DD}\ee
for each real positive value of the parameter $s$. Moreover, the generalized function $D$ so defined satisfies the association relations (\ref{DD}).
To show this, denote by
\be D_{\si}(x) \ := \ \frac{1}{\si} D\!\left(\si, \frac{x}{\si}\right), \  \h{where} \  \si = s(\pe). \label{strict}\ee
Now, for an arbitrary test-function $\ps\in\D$, \,evaluate the values
\[ I_1(\si) = \lan D_{\si}(x), \ps(x) \ran, \qquad  I_2(\si) = \lan D_{\si}^2 (x), \ps(x) \ran, \]
as $\ve \rar 0_+$, or equivalently, as $\si \rar 0_+$. But in view of (\ref{D-DD}), it is immediate to see that \ $\ \lim_{{}\si \rar 0_+} I_1(\si)
= \ \lim_{{}\si \rar 0_+} I_2(\si) = \lan \de, \ps \ran$; which according to Definition~2\,(b) implies (\ref{DD}).

\textbf{Remark.} The first equation in (\ref{DD}) is in consistency with the observation that $ D_{\si}(x) $ \,is a strict \de-net \,as defined in
distribution theory \cite[\S 7]{mob}. But note that $D$ is not the canonical embedding $\w{\de}$ of the \de-function since ${\w{\de}}^{\,2}$ does
not admit associated distribution.

The flexible approach to modelling singularities allowed by the generalized functions so that the models satisfy auxiliary conditions \,can be
systematically applied to defining generalized models of particular singularities. We recall that such approach has proved useful in studying
Euler-Lagrange equations for classical particle in \de-type potential as well as the geodesic equation for impulsive gravitational waves; see
\cite[\S 1.5, \,\S 5.3]{gct}.

We will consider models of singularities given by distributions with singular point support. For their definition, we intend to take advantage of
the properties of \de-modelling function $D$. Observe that it holds
\[ \vspace*{-2mm} (\de*D(s, .))(x)  = \lan \de_y, \,D(s, x-y) \ran  =  D(s, x), \]
\[ (\de'* D(s, .))(x)  = \lan \de'_y, \,D(s, x-y) \ran   = - \lan \de_y\, \,\pd_y D(s, x-y) \ran  = \lan \de_y, \,D'(s, x-y) \ran  =  D'(s, x). \]

\no This can be continued by induction for any derivative to define a generalized function $D^{(p)}(x)$ that models the distribution $\de^{(p)}(x)$
 and has a  representative
\be D^{(p)}(s, x) = (\de^{(p)}* D(s, .))(x). \label{D'}\ee

\no Clearly this is in consistency with the differentiation: \ $\pdx D^{(p)}(x) = D^{(p+1)}(x)$. Moreover,
\be D^{(p)}(-x) = \pdp D(-x) = (-1)^p D^{(p)}(x). \label{D-}\ee

In \cite{dam} we have employed such procedure for a unified modelling of singularities  given by distributions with singular point support, i.e.
 (besides $\de^{(p)}$) the distributions $ \nu_{\pm}^a, \,a\in\Om$. Namely, choosing an arbitrary generalized function $D$ with representative
 $D(s,x)$ that satisfies (\ref{D-DD}) for each $\vp \in \A$, we have introduced generalized functions $ X_{\pm}^{a}(x) $, modelling the above
 singularities, with representatives
\[ X_{\pm}^{a}(s, x) :=  (\nu_{\pm}^{a}*D(s, .))(x), \quad a\in\Om.\]
A consistency with the differentiation holds: \ $\pdx X_{\pm}^{a}(x)  =  X_{\pm}^{a-1}(x), \,\h{in particular,}\,H'  = D$, \ where $H\in\G$ is
model of the step-function $\te$, with representative $H(s,x) = \te *D(s, .)(x)$.

\vspace*{1mm}
Now we will define and study models in \G \ of functions on the real line that have discontinuities in a point. Their generalized models will be
obtained by the next definition that follows the idea suggested in \cite{dam}.

\vspace*{1mm}
\textsc{Definition}\,3. For any function $f \in C_d^k (\R\setminus\!\{0\})$, choosing an arbitrary function $D \in \G$ with representative $D(s,x)$
that satisfies (\ref{D-DD}) for each $\vp \in \A$, define its generalized model as the function $F(x) \in \G$ \ with representatives given by
\be F(s, x) :=  f*D(s, .)(x).\label{def}\ee

Note that each generalized functions $F$ so introduced is really model of the corresponding discontinuous function $f$. Indeed, let $F_{\si}(x)
= f*D_{\si}(x)$ be the representative of $F$ depending on $\si = \ve s$ and suppose (without loss of generality) that \ supp~$\!D(\si,x) \subseteq [-l, l]$ \,for some $l\in \R_+$. Then for an arbitrary test-function $\ps\in\D$, evaluate \  $I(\si) := \int_{\R}\ps(x) F_{\si}(x)\,dx.$
Transformation of the variable $y=\si v + x $ yields
\bea I(\si)& = & \frac{1}{\si}\int_{\R} dx\,\ps(x) \left[\int_{- \si l +x}^0 f(y) D\left(\si, \frac{x-y}{\si}\right)\,dy  + \int_0^{\,\si l +x }
f(y) D\left(\si, \frac{x-y}{\si}\right)\,dy \right] \\ & = & \int_{\R} dx\,\ps(x) \left[\int_{-l}^{-x/\si} f(\si v+x) D(\si, v)\,dv  +
\int_{-x/\si}^{\,l} f(\si v+x) D(\si, v)\,dv \right] \eea
Now taking the limit as $\ve \rar 0_+$, or else $\si =\ve s \rar 0_+ $, and applying equation (\ref{D-DD}), we get
\[ \lim_{\si \rar 0_+} I(\si) \ = \  \int_{\R} dx\,\ps(x)\,f(x) \int_{-l}^{\,l} D(\si, v)\,dv = \lan\,f, \ps\, \ran,\]
which according to Definition~2(b) implies the association \  $F(x) \approx f(x)$.

\vspace*{5mm}
\no \textbf{4 \ PRODUCTS OF SOME SINGULARITIES MODELLED IN} {\boldmath{\GR}\unboldmath}

\vspace*{2mm}
The models of singularities we consider all have products in Colombeau algebra as generalized functions, but we are seeking results that can be
evaluated back in terms of distributions, i.e. products that admit associated distributions.

So, it was proved in \cite{dam}\,that, for an arbitrary $p$ in $\N_0$, the generalized models $X_{\pm}^{\,p},  D^{(p+1)}$, and $D^{(p+2)}$ satisfy:
\be
(\mp 1)^{p+1} X_{\pm}^{{}p}\,.\,D^{(p+1)}\ \approx \ \de \mp \frac{(p+1)}{2}\de'. \label{d+1}
\ee
\be
(\mp 1)^{p}\, X_{\pm}^{{}p}\,.\,D^{(p+2)} \ \approx \mp \frac{2p+3}{2} \de'  + \frac{1}{2}  \left(\begin{array}{c}\!\!p+2\!\!\\ \!\!2\!\!\end{array}
\right) \de'', \label{d+2} \ee

\no Then the generalized functions $ |X|^{\,p} = X_{+}^{\,p} + X_{-}^{\,p}, \ |X|^{\,p}\,\mathrm{sgn}\,x  = X_{+}^{\,p} -  X_{-}^{\,p}$, and $D$
satisfy:
\be |X|^{\,p}\,.\,D^{(p+1)} \approx  2\,\de, \quad  |X|^{\,p}\,.\,D^{(p+2)} \ \approx \  (2p+3)\,\de', \ \h{for} \  p= 1, 3, 5,... \label{odd} \ee
\be |X|^{\,p}\,\mathrm{sgn}\,x \,.\,D^{(p+1)} \ \approx \  - 2\,\de, \quad  |X|^{\,p}\,\mathrm{sgn}\,x \,.\,D^{(p+2)} \ \approx \ - (2p+3)\,\de', \
\h{for} \ p = 2, 4, 6... . \label{even} \ee

\vspace*{2mm} Now we proceed to studying singular products,
obtained by Definition 3, of generalized models of functions from
the class $ C_d^k(\R\setminus\!\{0\})$ with derivatives of the
\de-modelling generalized function $D(x)$. Products of
discontinuous functions with the derivatives of \de \ exist
neither in the classical Distribution theory nor as so-called
Colombeau products - their canonical embeddings in \G \ do not
admit associated distributions. Nevertheless, their generalized
models obey the following.

\no \begin{Th1}$\!\!.$   For each function $f(x) \in C_d^2(\R\setminus\!\{0\})$,  its model $F(x)$ in \GR \ satisfies
\be F(x)\,.\,D'(x) \ \approx \  - ( h_0 + m_1)\ \de \ + \ m_0 \ \de',\label{fd'}\ee
where $ h_0 = f(0_+) - f(0_-), m_0 = ( f(0_-) + f(0_+) ) / 2, \,m_1 = ( f'(0_-) + f'(0_+) ) / 2.$
\end{Th1}

\no \textit{Proof:} For $\ps(x) \in \D$, we denote \ $I(\si) := \lan\  F_{\si}(x)\,.\,D'_{\si} (x), \,\ps(x)\,\ran $. From equations
(\ref{strict}), (\ref{D'}), and (\ref{def}),  we get on transforming the variables  $y= \si v + x, \,x= - \si u$ and taking into account
equation (\ref{D-})
\[ I(\si) =  \frac{- 1}{\si} \int_{- l}^{\,l} du\,\ps(- \si u)\,D'(\si, u) \left(\int_{-l}^{\,u} f(\si v -\si u) D(\si, v)\,dv  + \int_{\,u}^{\,l}
 f(\si v -\si u) D(\si, v)\,dv \right) \]
Now applying Taylor theorem  to the test-function \ps \ and changing the order of integration, we obtain
\bean I(\si) & = & \frac{- \,\ps(0)}{\si} \int_{- l}^{\,l} dv\, D(\si, v) \left(\int_v^{\,l} f(\si v -\si u) D'(\si, u)\,du + \int_{- l}^{\,v}
f(\si v -\si u) D'(\si, u)\,du\right) \nonumber \\ & &
\!\!\! + \,\ps'(0) \!\int_{- l}^{\,l} dv\,D(\si, v) \!\left(\int_v^{\,l} f(\si v -\si u) u D'(\si, u)\,du  + \!\int_{- l}^{\,v} f(\si v -\si u) u
D'(\si, u)\,du\right) \nonumber \\*[2mm] & &
\!\!\! + \ O(\si) \ =: \ \ps(0) \ I_1\ + \ \ps'(0) \ I_2 + \ O(\si). \label{11}\eean
To obtain the above asymptotic evaluation, we have taken into account that the third term in the Taylor expansion is multiplied by definite
integrals majorizable by constants.

\no Integrating further by parts in the variable $u$, applying Lebesque theorem on bounded convergence, taking into account equation (\ref{D-DD}),
and making use of Notation~3, we obtain
\bea I_1 & = & \frac{f(0_-)}{\si} \int_{- l}^{\,l} D^2(\si, v)\,dv -
\int_{- l}^{\,l} dv\, D(\si, v) \int_v^{\,l} f'(\si v -\si u) D(\si, u)\,du \\ & & -  \frac{f(0_+)}{\si} \int_{- l}^{\,l} D^2(\si, v)\,dv -
\int_{- l}^{\,l} dv\, D(\si, v) \int_{-l}^{\,v} f'(\si v -\si u) D(\si, u)\,du \\ & = &
\,f(0_-) - f(0_+) - \frac{1}{2}\,\left[\,f'(0_-) + f'(0_+)\,\right]  + \ O(\si) \ = \ - \,(h_0 + m_1) \ + \ O(\si).\eea
Here the Taylor theorem is applied to the function $ f' \in C^1(\R_{\pm})$ up to second order to get its expansion about the point $\si (v -u)$
which is respectively $>0$ or $< 0$. This gives
\[ f'(\si v -\si u) = f'(0_{\pm}) + \si (v -u)\ f''(0_{\pm}) + \ O(\si^2) = f'(0_{\pm}) + \ O(\si), \ \h{for} v > u, \h{resp.,} \,v < u.\]
Proceeding similarly as above, we get for the second term in (\ref{11})
\bea I_2 & = & - \,f(0_-) \int_{- l}^{\,l} v D^2(\si, v)\,dv - f(0_-)
\,\int_{- l}^{\,l} dv\, D(\si, v) \int_v^{\,l} D(\si, u)\,du \\ & & + \,f(0_+) \int_{- l}^{\,l} v D^2(\si, v)\,dv - f(0_+)\,\int_{- l}^{\,l} dv\,
D(\si, v) \int_{-l}^{\,v} D(\si, u)\,du + \ O(\si) \\ & = &
 \frac{- 1}{2}\,\left[\,f(0_-) + f'(0_+)\,\right]  + \ O(\si) \ = \ - \,m_0 \ + \ O(\si). \eea
We have used that $D(\si, v)$ satisfies
\[ \int_{- l}^{\,l} v D^2(\si, v)\,dv = 0 \ \h{and} \ \int_{- l}^{\,l} dv\,D(\si, v) \int_{- l}^{\,v}D(\si, u)\,du = \frac{1}{2}. \]
Therefore
\bea I(\si) & = & \lan\  F_{\si}(x)\,.\,D'_{\si} (x), \,\ps(x)\,\ran =  - \,(h_0 + m_1) \,\ps(0) \ - \ m_0 \,\ps'(0) \ + \ O(\si) \\ & = & -
\,(h_0 + m_1) \,\lan\,\de, \ps\,\ran \ + \ m_0 \,\lan\,\de', \ps\,\ran \ + \ O(\si). \eea
Passing then to the limit as  $\si  \rar 0_+ $ and applying Definition~2\,(b), we obtain equation (\ref{fd'}).

\vspace*{2mm}
Next we consider the product of models of discontinuous functions with the second derivative of the \de-modelling function $D(x)$. In view of
Notation~3, one has the following.

\vspace*{1mm}
\no \begin{Th2}$\!\!.$   For each function $f(x) \in C_d^3(\R\setminus\!\{0\})$,  its generalized model $F(x)$ \,in \,\GR \  satisfies\,:
\be F(x)\,.\,D''(x) \ \approx \  (h_1 + m_2)\,\de - \left( \frac{3}{2}\,h_0 + 2 m_1\right)\,\de' + m_0\,\de''.
\label{fd''}\ee
\end{Th2}

\vspace*{1mm}
\no \textit{Proof:} For $\ps(x) \in \D$, we denote \ $J(\si) := \lan\  F_{\si}(x)\,.\,D''_{\si} (x), \,\ps(x)\,\ran $. From equations (\ref{strict}),
(\ref{D'}), and (\ref{def}),  we get on transforming the variables  $y= \si v + x, \,x= - \si u$ and taking into account equation (\ref{D-})
\[J(\si) = \!\frac{1}{\si^2} \!\int_{- l}^{\,l} du\,\ps(- \si u)\,D''(\si, u) \!\left(\!\int_{-l}^{\,u} f(\si v -\si u) D(\si, v)\,dv  +
\!\!\int_{\,u}^{\,l} f(\si v -\si u) D(\si, v)\,dv \right) \]
Now applying Taylor theorem to the test-function \ps \ and changing the order of integration, we obtain
\bean J(\si) & \!\! = & \!\!\!\frac{\ps(0)}{\si^2} \!\int_{- l}^{\,l} dv\,D(\si, v) \left(\int_v^{\,l} f(\si v -\si u) D''(\si, u)\,du +
\!\int_{- l}^{\,v} \!f(\si v -\si u) D''(\si, u)\,du\right) \nonumber \\ & &
\!\!\! - \frac{\ps'(0)}{\si} \!\!\int_{- l}^{\,l} \!dv\,D(\si, v) \!\left(\int_v^{\,l} f(\si v -\si u) u D''(\si, u)\,du  \!+
\!\int_{- l}^{\,v}\!f(\si v -\si u) u D''(\si, u)\,du\right) \nonumber \\*[2mm] & & \!\!\! + \frac{\ps''(0)}{2} \!\!\int_{- l}^{\,l}\!dv\,D(\si, v)
\!\!\left(\int_v^{\,l}\!f(\si v -\si u) u^2 D''(\si, u)\,du  + \!\!\int_{- l}^{\,v}\!f(\si v -\si u) u^2 D''(\si, u)\,du\right) \nonumber \\*[2mm] & &
\!\!\! + \ O(\si) \ =: \ \ps(0) \ J_1\ - \ \ps'(0) \ J_2 + \ \frac{\ps''(0)}{2} \ J_3 \ + \ O(\si). \label{22}\eean

\no Integrating twice by parts in the variable $u$, applying Lebesque theorem on bounded convergence, and taking account of equation (\ref{D-DD}),
we get further
\bea J_1 & = & \frac{- f(0_-)}{\si^2} \int_{- l}^{\,l} D(\si, v) D'(\si, v)\,dv  + \frac{1}{\si} \int_{- l}^{\,l} dv\, D(\si, v) \int_v^{\,l}
f'(\si v -\si u) D'(\si, u)\,du \\ & & + \frac{f(0_+)}{\si^2} \int_{- l}^{\,l} D(\si, v) D'(\si, v)\,dv  + \frac{1}{\si}\int_{- l}^{\,l} dv\,D(\si, v)
 \int_{- l}^{\,v}f'(\si v -\si u) D'(\si, u)\,du \\ & = &
- f'(0_-) \frac{1}{\si} \int_{- l}^{\,l} D^2(\si, v) \,dv +  \int_{- l}^{\,l} dv\, D(\si, v) \int_v^{\,l} f''(\si v -\si u) D(\si, u)\,du  \\ & & +
 f'(0_+) \frac{1}{\si} \int_{- l}^{\,l} D^2(\si, v) \,dv + \int_{- l}^{\,l} dv\, D(\si, v) \int_{- l}^{\,v} f''(\si v -\si u) D(\si, u)\,du  + \ O(\si)
 \\ & = &
- \,f'(0_-) + f'(0_+) + \,\frac{1}{2}\,\left[\,f''(0_-) + f''(0_+)\,\right]  + \ O(\si) \ = \  \,(h_1+ m_2) \ + \ O(\si).\eea

\no Here the Taylor theorem is applied to the function $ f'' \in C^1(\R_{\pm})$ up to second order to get its expansion about the point $\si (v -u)$
which is respectively $>0$ or $> 0$, obtaining thus
\[ f''(\si v -\si u) = f''(0_{\pm}) + \ O(\si).\]
Throughout the calculations, we have used that, due to $ D(., -x ) = D(., x)$, it holds\,:
\[ \int_{- l}^{\,l} D(\si, v) D'(\si, v)\,dv = 0 \ \h{and} \ \int_{- l}^{\,l} dv\,D(\si, v) \int_{- l}^{\,v}D(\si, u)\,du = \frac{1}{2}. \]
Proceeding similarly as above, we get for the second and third terms in (\ref{22})\,:
\bea J_2 & = & \frac{- \,f(0_-)}{\si} \int_{- l}^{\,l} v D(\si, v) D'(\si, v)\,dv
+ \frac{f(0_+)}{\si} \int_{- l}^{\,l} v D(\si, v) D'(\si, v)\,dv \\ & &
- \frac{1}{\si}\int_{- l}^{\,l}dv\,D(\si, v) \left(\int_{\,v}^{\,l} f(\si v -\si u) D'(\si, u)\,du
+  \int_{- l}^{\,v} f(\si v -\si u) D'(\si, u)\,du \right) \\ &  &
+ \int_{- l}^{\,l} dv\,D(\si, v) \left( \int_v^{\,l}f'(\si v -\si u) u D'(\si, u)\,du
+  \int_{- l}^{\,v} f'(\si v -\si u) u D'(\si, u)\,du \right) \\ & = &
\frac{1}{2}\,f(0_-) - \frac{1}{2}\,f(0_+) + f(0_-) \frac{1}{\si}\int_{- l}^{\,l}D^2(\si, v)\,dv - f(0_+) \frac{1}{\si}\int_{- l}^{\,l}D^2(\si, v)
\,dv \\ & &
- 2 \,f'(0_-) \int_{- l}^{\,l} dv\,D(\si, v) \int_v^{\,l} D(\si, u)\,du - 2 \,f'(0_+) \int_{- l}^{\,l} dv\,D(\si, v) \int_{- l}^{\,v} D(\si, u)\,du
\\ & &
\!\!\!+ 2 \!\int_{- l}^{\,l} dv\,D(\si, v) \!\left(\int_v^{\,l} f(\si v -\si u) u D'(\si, u)\,du  + 2 \!\int_{- l}^{\,v} f(\si v -\si u) u D'(\si, u)
\,du\right) \\ & &
\\ & = &
 \frac{3}{2} f(0_-) - \frac{3}{2} f(0_+)  - \,f'(0_-) - \,f'(0_+) + \ O(\si) = - \left(\,\frac{3}{2} h_0  + 2 \,m_1\right) \ + \ O(\si).\eea

\bea J_3 & = & - \,2 \int_{- l}^{\,l} dv\,D(\si, v) \!\left(\int_v^{\,l} f(\si v -\si u) u D'(\si, u)\,du  +  \int_{- l}^{\,v} f(\si v -\si u) u
 D'(\si, u)\,du\right) \\ & &
- f(0_-) \int_{- l}^{\,l} v^2 D(\si, v) D'(\si, v)\,dv + f(0_+) \int_{- l}^{\,l}v^2 D(\si, v) D'(\si, v)\,dv  + \ O(\si) \\ & = &
 \!\!\!2\,f(0_-) \!\int_{- l}^{\,l} dv\,D(\si, v) \int_v^{\,l} D(\si, u)\,du  + f(0_+)  \!\int_{- l}^{\,l} dv\,D(\si, v) \int_{- l}^{\,v} D(\si, u)
 \,du +  O(\si)\\ & = &
f(0_-) \ + \ f(0_+) \ + \ O(\si) \ = \ 2\,m_0 \ + \ O(\si).\eea

\no During the calculations of $J_2$ and $J_3$, we have used that
\[ \int_{- l}^{\,l} v D^2(\si, v) \,dv = \int_{- l}^{\,l} v^2 D(\si, v) D'(\si, v)\,dv = 0 \ \h{and} \ \frac{1}{\si} \int_{- l}^{\,l} v D(\si, v)
D'(\si, v)\,dv = - \frac{1}{2}.\]

\no Thus we finally get
\bea \lan F_{\si}(x)\,.\,D''_{\si} (x), \,\ps(x)\ran & = & \!\!\!(h_1 + m_2)\,\ps(0) + \left( \frac{3}{2}\,h_0 + 2 m_1\right)\,\ps'(0) + m_0\,\ps''(0)
 + O(\si) \\ & = & \!\!\!(h_1 + m_2) \,\lan\,\de, \ps\,\ran - \left( \frac{3}{2}\,h_0 + 2 m_1\right) \,\lan\,\de', \ps\,\ran + m_0\,\lan\,\de'',
  \ps\,\ran + O(\si). \eea

\no Then passing to the limit as  $\si  \rar 0_+ $  and applying Definition~2\,(a), we obtain equation (\ref{fd''}).

\vspace*{1mm}
\no \textbf{Remarks. (1)} \,In the particular case of infinitely-differentiable functions, equations (\ref{fd'}) and (\ref{fd''}) are in consistence
with the corresponding products obtained in classical Distribution theory.  Indeed, for $ f \in C^{\infty}(\R)$ and $\ps \in \D$, we have
\[ \lan\,f\,.\,\de', \,\ps\,\ran =  \lan\,\de', \, f\,\ps\,\ran = - \pdx\,(\,f\,\ps)\left|_{x=0}\right.,  \ \h{or else,} \quad  f\,.\,\de' \ =
 \ -f'(0)\,\de + f(0)\,\de'.\]
This clearly coincides with equation (\ref{fd'}) since in this case $ m_0 = f(0), m_1 = f'(0), h_0 =0$. Similar argument applies to equation
(\ref{fd''}).

\textbf{(2)} In the results obtained so far, we have restricted ourselves to studying singularities in the point $x=0$, but the considerations
 are clearly valid for any other singular point.

\vspace*{1mm}
Recall next that any function can be canonically represented as a sum of its even and odd parts\,:
\[f(x) = \sum_{\sigma = 0,1}f_{\sigma}(x), \quad \h{where} f_0(x) := \frac{1}{2}\ [ f(x) + \check{f}(x)] \ \h{and} f_1(x) := \frac{1}{2}\ [ f(x) -
\check{f}(x)] \]
\no are indeed even and odd functions\,:  \ $f_{\sigma}(-x) = (-1)^ {\sigma} \,f_{\sigma}(x), \ \sigma = (0,1).$

Then, as consequences from the results of Theorems~1--2, we obtain this.

\no \begin{Cor1}$\!\!.$  The even and odd parts $F_0(x), F_1(x)$ of the generalized model $F(x) $ \,of each function $f(x) \in C_d^3
(\R\setminus\!\{0\})$ \,satisfy the equations\,:
\be F_0(x)\,.\,D'(x) \ \approx \  m_0 \ \de'(x), \qquad \qquad F_1(x)\,.\,D'(x) \ \approx \  - ( h_0 + m_1)\ \de(x). \label{eofd'}\ee
\be F_0(x)\,.\,D''(x) \ \approx \  (h_1 + m_2)\,\de  + m_0\,\de'', \qquad F_1(x)\,.\,D''(x) \ \approx \ - \left( \frac{3}{2}h_0 + 2 m_1\right)\,\de'.
\label{eofd''}\ee
\end{Cor1}

\no \textit{Proof:}\ Rewrite equations (\ref{fd'}) and (\ref{fd''})  for the even and odd parts of the model $F(x)$. Then the above equations follow
on the observation that the derivative of an even/odd function is, respectively, odd/even function, and also that $h(f_0) = m(f_1) = 0$.

\newpage
\textbf{Examples.} The following equations can be obtained on replacing the generalized model $F(x)$ of the function $f(x)\in C_d^3
(\R\setminus\!\{0\})$ successively with\,:

\textbf{(a)} \ $H $ \ in the first equation of (\ref{fd'})\,: $ \quad \qquad \qquad\qquad \qquad H\,.\,D' \ \approx \ - \,\de + \frac{1}{2}\ \de'. $

\textbf{(b)} \ $H$ \ in the first equation of (\ref{fd''})\,: $ \quad \qquad\qquad \qquad \qquad H\,.\,D'' \ \approx \ - \, \frac{3}{2}\ \de' +
\frac{1}{2}\ \de''.$

\textbf{(c)} \ $X_+$  \ in the first equation of (\ref{fd''})\,: $ \qquad \qquad \qquad \qquad X_+\,.\,D'' \ \approx \ \de$ - \de'.

\textbf{(d)} \ $|X|$ \ in the first equation of (\ref{eofd''})\,:  $ \qquad \qquad \qquad \qquad |X|\,.\,D'' \ \approx  \ 2\,\de$.

\textbf{(e)} \ $|X|\,\mathrm{sgn}\,x$  \ in the second equation of (\ref{eofd''})\,:  $ \qquad |X|\mathrm{sgn}\,x \,.\,D'' \ \approx \ - \,2\,\de$.

\vspace*{2mm}
\no Note that these equations coincide correspondingly with\,: \, \textbf{(a)} \,equation (\ref{d+1}) for $p=0$, \,\textbf{(b)}~\,equation
(\ref{d+2}) for $p=0$, \,\textbf{(c)} \,equation (\ref{d+1}) for $p=1$,  \,\textbf{(d)} \,the first equation in (\ref{odd}) for $p=1$,
\,\textbf{(e)} \,the first equation in (\ref{even}) for $p=1$.

\vspace*{1mm}
\setlength{\baselineskip}{15pt}


\begin{thebibliography}{19}

\bibitem{col84} Colombeau, J.-F. (1984) \textit{New Generalized Functions and Multiplication of Distributions}, North Holland Math. Studies 84,
 Amsterdam.

\bibitem{col87} Colombeau, J.-F., A.Y. Le Roux. (1987) \textit{ In: Advances in Comp. Methods for PDEs~VI}, R.Vichnevetsky, R.S. Steplemen (editors),
Publ. IMACS, 28--37.

\bibitem {col89} Colombeau, J.-F., A.Y. Le Roux, A.Nuissair, B.Perrot. (1989) Microscopic profiles of shock waves and ambiguities in multiplication
 of distributions. \textit{SIAM J.Num.Anal}. \textbf{26}, 87--883.

\bibitem{dam} Damyanov, B. Moddeling and products of singularities in Colombeau algebra, J. Applied Analysis,
14 (2008), 1, 89--102.

\bibitem{gct} Grosser, M., M. Kunzinger, M. Oberguggenberger, R. Steinbauer. \,\textit{Geometric Theory of Generalized Functions with
Applications to General Relativity}, Springer Series Mathematics
and Its Applications , Vol. 537, 2002.

\bibitem{HO} Hermann, R., M. Oberguggenberger. (1999) ODEs and generalized functions. \textit{ In: Nonlinear Theory of Generalized Functions}.
M.Grosser, G.H\"ormann, M.Kunzinger, M. Oberguggenberger (editors), Chapman \& Hall/CRC, Boca Raton.

\bibitem{mob} Oberguggenberger M. (1992) \textit{Multiplication of Distributions and Applications to PDEs}, Longman, Essex.

\bibitem{St} Steinbauer, R. (1997) Geodesics and geodesic derivation for impulsive waves. \textit{J. Math. Phys.} \textbf{38}, 1614--1622.

\end{thebibliography}
\end{document}